\documentclass[article,11pt,reqno]{amsart}
\usepackage[reqno]{amsmath}
\usepackage[utf8]{inputenc}
\usepackage{graphicx, amsmath, amssymb, amscd, amsthm, euscript, psfrag, amsfonts,bm}
\usepackage{enumitem}
\usepackage{amsmath}
\usepackage{etoolbox}
\usepackage{amssymb}
\usepackage{amscd}
\usepackage{amsthm}
\usepackage{amsfonts}
\usepackage{graphicx}
\usepackage{tensor}
\usepackage{multirow}
\usepackage{multicol}
\usepackage{mathscinet}
\usepackage[T1]{fontenc}
\usepackage{mathrsfs}
\usepackage{appendix}
\usepackage{tikz}
\usetikzlibrary{calc,math,intersections}
\usepackage{float}

\usepackage[colorlinks=true]{hyperref}
\hypersetup{urlcolor=blue, citecolor=red}
\usepackage{cleveref}
\crefformat{section}{\S#2#1#3} 
\crefformat{subsection}{\S#2#1#3}
\crefformat{subsubsection}{\S#2#1#3}

\newtheorem{theorem}{Theorem}[section]
\newtheorem*{theorem*}{Theorem}
\newtheorem{corollary}[theorem]{Corollary}

\newtheorem{lemma}[theorem]{Lemma}

\newtheorem{proposition}[theorem]{Proposition}

\newtheorem*{problem}{Problem}
\theoremstyle{definition}

\newtheorem{remark}[theorem]{Remark}

\title{ On identities concerning integer parts}

\author{Zichang Wang}
\address{Qiuzhen College, Tsinghua University, Beijing, 100084, China}
\email{w-zc21@mails.tsinghua.edu.cn}
\date{}
\author{Chengyang Wu}
\address{School of Mathematical Sciences, Peking University, Beijing, 100871, China
}
\email{2101110022@stu.pku.edu.cn }
\author{Bohan Yang}
\address{Yau Mathematical Sciences Center, Tsinghua University, Beijing, 100084, China }
\email{ybh20@mails.tsinghua.edu.cn}

\begin{document}

\begin{abstract}
 In 2007 V. Zhuravlev discovered a family of identities concerning integer parts which are satisfied by the number
 $\frac{\sqrt{5}+1}{2}$. Some of these identities turned out to be characterization properties
 of the number  $\frac{\sqrt{5}+1}{2}$.  
 In this paper we generalize the simplest of these identities.

\end{abstract}

\maketitle
\section{Introduction}

   For any real number $\alpha$, let $[\alpha]$ denote its integer part, and $\{\alpha\}=\alpha-[\alpha]$ denote its fractional part. By studying the first recurrence map, V. Zhuravlev \cite{z}  discovered a family of identities concerning integer parts which are satisfied by the number $\alpha=\frac{\sqrt{5}+1}{2}$. The two simplest identities in this family are
\begin{equation}\label{z1}
[([n\alpha]+1)\alpha] = [n\alpha^2]+1,\,\,\, n \in \mathbb{Z}
\end{equation}
and 
\begin{equation}\label{z2}
[[n\alpha]\alpha]+1 = [n\alpha^2],\,\,\,\, n \in \mathbb{Z}\setminus\{0\}.
\end{equation}
It should be noted that, for $n\in\mathbb{Z}\setminus\{0\}$ the substitution $ n \mapsto -n$ changes (\ref{z1}) to (\ref{z2})
and vice versa. One can find in \cite{av} more general statements about rotations of a circle which are related to such identities concerning integer parts.
Later A. Zhukova and A. Shutov in \cite{S21} proved that identities  (\ref{z1}) and (\ref{z2}) actually characterize the number $\alpha$. Namely, they proved the following result:
\begin{theorem}\label{0000}
        A real number $\alpha$ satisfies the identity
    \begin{equation}\label{idfor2}
     \left[ \left[  n\alpha \right]\alpha \right]+1 =  \left[   n\alpha^2 \right] 
    \end{equation}
    for all $n\in\mathbb{Z}\setminus \{0\}$, if and only if $\alpha=\frac{\sqrt{5}+1}{2}$.
\end{theorem}

In fact, A. Zhukova and A. Shutov in \cite{S21} proved a stronger statement. They showed that to characterize the real number $\alpha =\frac{\sqrt{5}+1}{2}$ it is sufficient to show that the identity (\ref{idfor2}) holds for all $n$ in the Fibonacci sequence.
General and specific problems related to combinatorics of sequences about integer parts have been considered in the fundamental paper \cite{aa}.

\vskip+0.3cm
In this paper we generalize Theorem \ref{0000} and consider some other identities which also characterize certain algebraic integers.

\vskip+0.3cm

Let us start with the following result which contains Theorem \ref{0000} as a special case when $l=k=1$. We will give it a short dynamical proof using equidistribution in Section \ref{proofofmain}.

\begin{theorem} \label{main}
Let $l$ and $k$ be two positive integers. 
Then a positive real number $\alpha$ satisfies 
\begin{equation}\label{algforkl}
    \alpha^{l+k}-\alpha^{l}\in \mathbb{Z}\cap [1,2^{\frac{l}{k}})
\end{equation}
if and only if the identity
\begin{equation}\label{idforkl}
    [[n\alpha^l]\alpha^k]+1 = [n\alpha^{l+k}]
\end{equation}
holds for all $n\in\mathbb{Z}\setminus\{0\}$.
\end{theorem}
\vskip+0.3cm

We now formulate three partial generalizations of Theorem \ref{main}. Their proofs are similar to the one of Theorem \ref{main}. The first of them replaces $+1$ on the left-hand side of (\ref{idforkl}) by $+\delta$ for some $\delta$ slightly smaller than $1$.

\begin{theorem} \label{gen1}
Let $l$ and $k$ be two positive integers, and let $\alpha\in (0,1)\cup (1,2^{\frac{1}{k}})$
be a real number.
Then the following statements are equivalent:
\begin{enumerate}
    \item[\rm(i)] $\alpha$ satisfies 
\begin{equation*}
    \alpha^{l+k}-\alpha^{l}\in \mathbb{Z}\cap [1,2^{\frac{l}{k}}).
\end{equation*} 
\item[\rm(ii)] 
there exists some $\delta\in [\alpha^k-1,1)$ such that
the identity
\begin{equation*}
    [[n\alpha^l]\alpha^k+\delta] = [n\alpha^{l+k}]
\end{equation*}
holds for all $n\in\mathbb{Z}$.
\item[\rm(iii)] for any $\delta\in [\alpha^k-1,1)$,
the identity
\begin{equation*}
    [[n\alpha^l]\alpha^k+\delta] = [n\alpha^{l+k}]
\end{equation*}
holds for all $n\in\mathbb{Z}$.
\end{enumerate}
\end{theorem}

\begin{remark}
    It is trivial to see that $\alpha=1$ satisfies (ii) and (iii) but not (i) in Theorem \ref{gen1}. This will be reflected in the proof of Theorem \ref{gen1}. Moreover, it is also trivial to check that the condition $\delta<1$ makes the identity also true for $n=0$, which is not allowed in Theorem \ref{main}.
\end{remark}
 
 The second partial generalization of Theorem \ref{main} replaces the coefficient $-1$ of $\alpha^l$ on the left-hand side of (\ref{algforkl}) by any negative integer $-m$.

\begin{theorem}\label{gen2}
Let $l,k,$ and $m$ be positive integers. Then a positive real number $\alpha$ satisfies
\begin{equation}\label{algforllm}
    \alpha^{l+k}-m\alpha^l\in \mathbb{Z}\cap [1,(m+1)^{\frac{l}{k}}),
\end{equation}
if and only if the identity
    \begin{equation}\label{idforllm}
    [[n\alpha^l]\alpha^k]+[nm\alpha^l]+1-m[n\alpha^l]=[n\alpha^{l+k}]
\end{equation}
holds for all $n\in\mathbb{Z}\setminus\{0\}$. 
\end{theorem}

The third partial generalization of Theorem \ref{main} restricts the identity (\ref{idforkl}) for a particular subsequence $\{a(n)\}_{n\in\mathbb{Z}}$ of nonzero integers, where one of its key properties is that the orbit $\{a(n)\gamma \pmod{\mathbb{Z}}\}_{n\in\mathbb{Z}}$ is dense in $\mathbb{T}$ for all irrational numbers $\gamma$. However, we can say nothing for a general subsequence without this property, say for the Fibonacci sequence mentioned in \cite{S21}.

\begin{theorem} \label{gen3}
Let $l$ and $k$ be two positive integers, $P(X)\in\mathbb{Z}[X]$ be a non-constant polynomial,
and $\alpha$ be a positive real number with
$\{\alpha^l,\alpha^{k}\}\not\subseteq\mathbb{Q}$.
Then $\alpha$ satisfies 
\begin{equation*}
    \alpha^{l+k}-\alpha^{l}\in \mathbb{Z}\cap [1,2^{\frac{l}{k}})
\end{equation*}
if and only if the identity
\begin{equation}\label{idforpoly}
    [[P(n)\alpha^l]\alpha^k]+1 = [P(n)\alpha^{l+k}]
\end{equation}
holds for all $n\in\mathbb{Z}\setminus\{
\text{roots of } P(X)
\}$.
\end{theorem}

\begin{remark}
    The additional assumption $\{\alpha^l,\alpha^{k}\}\not\subseteq\mathbb{Q}$ is necessary. 
    Otherwise, let us take $l=k=1$, $P(X)=4X+1$, and $\alpha=\frac{3}{2}$.
    Then it is straightforward to see that $\alpha$ doesn't satisfy 
    $$\alpha^2-\alpha\in\mathbb{Z}\cap [1,2),$$
    but the identity
    $$
    [[P(n)\alpha]\alpha]+1 = [P(n)\alpha^{2}]
$$
holds for all $n\in\mathbb{Z}$. 
\end{remark}

It is natural to ask if there is an analogy of above theorems concerning three or more nested brackets. We find it much more complicated and raise the following problem:

\begin{problem} Is there any algebraic integer $\alpha$ characterized by the identity
    \begin{equation}\label{formula}
    [[[n\alpha]\alpha]\alpha]+1=[n\alpha^3]
\end{equation}
for all $n\in\mathbb{Z}\setminus\{0\}$?
\end{problem}

Our paper is organized as follows. In Section \ref{proofofmain} we give a short but complete proof of Theorem \ref{main}.
In Sections \ref{proofofgen1} - \ref{proofofgen3} we sketch the proofs of the three partial generalizations of Theorem \ref{main}. In the Appendix we state and prove the well-known Weyl's well-known equidistribution theorem for our needs.

\section{Proof of Theorem \ref{main}} \label{proofofmain}

Let us start with the easier direction of Theorem \ref{main} based on the following fact.

\begin{lemma}\label{samedenom}
    Let $\alpha$ be a real number with $\{\alpha^l,\alpha^k\}\subseteq\mathbb{Q}$ for two positive integers $l,k$. Write $\alpha^l=\frac{p}{q}$ and
    $\alpha^k=\frac{p'}{q'}$, where
    $p,p'\in\mathbb{Z}$, $q,q'\in\mathbb{N}$, and $\mathrm{gcd}(p,q)=\mathrm{gcd}(p',q')=1$. 
    Then $\mathrm{gcd}(q,q')\neq 1$ unless $q=q'=1$.
\end{lemma}
\begin{proof}
Write $d=\mathrm{gcd}(l,k)$. By applying the Euclidean algorithm to the pair $(l,k)$, we see by induction that $\alpha^d$ is of exactly one of the following forms:
$$
\frac{p^{a}(q')^{b}}{q^{a}(p')^{b}}
\text{ or }
\frac{q^{a}(p')^{b}}{p^{a}(q')^{b}}
$$
for some $a,b\in\mathbb{N}\cup\{0\}$.
Without loss of generality, we may assume $\alpha^d$ is of the latter form.
It follows that $\alpha^l$ is also of the latter form, namely, for $c=\frac{l}{d}\in\mathbb{N}$, we have
$$\alpha^l=\frac{q^{ac}(p')^{bc}}{p^{ac}(q')^{bc}}.$$
Then we have $q\mid p^{ac}(q')^{bc}$ since $\mathrm{gcd}(p,q)=1$,
and hence $q\mid (q')^{bc}$.

Suppose that $\mathrm{gcd}(q,q')=1$. Then it follows that $q=1$, namely, $\alpha^l\in\mathbb{Z}$. Since $\alpha^d$ is rational and $(\alpha^{d})^c=\alpha^l$ is an integer, we conclude that $\alpha^d$ must be an integer and hence $\alpha^k$ is also an integer, namely, $q'=1$.
\end{proof}

\begin{proof}[Proof of the ``only if'' part of Theorem \ref{main}]
It is easy to deduce that 
\begin{equation}\label{rangeofalpha}
1<\alpha<2^{\frac{1}{k}}
\end{equation}
from the conditions $\alpha>0$ and $1\leq \alpha^l(\alpha^k-1)<2^{\frac{l}{k}}$.
Moreover, let us write
\begin{equation}\label{algfort}
\alpha^{l}(\alpha^k-1)=M
\end{equation}
for some $M\in \mathbb{Z}\cap [1,2^{\frac{l}{k}})$ from (\ref{algforkl}).
We claim that 
$\{\alpha^l,\alpha^k\}\cap\mathbb{Q}=\varnothing$.

In fact, suppose to the contrary that $\{\alpha^l,\alpha^k\}\cap\mathbb{Q}\neq\varnothing$.
Then it follows from (\ref{algfort}) that 
$\{\alpha^l,\alpha^k\}\subseteq\mathbb{Q}$.
Let us write $\alpha^l=\frac{p}{q}$ and
    $\alpha^k=\frac{p'}{q'}$, where $\mathrm{gcd}(p,q)=\mathrm{gcd}(p',q')=1$.
Since $\alpha^k\in (1,2)$ cannot be an integer,
we have $\mathrm{gcd}(q,q')\neq 1$ by Lemma \ref{samedenom}.
However, it can be directly calculated from (\ref{algfort}) that
$$
\frac{p}{q}=\frac{Mq'}{p'-q'},
$$
which implies that $q\mid (p'-q')$.
In particular, this means that 
$\mathrm{gcd}(q,q')\mid \mathrm{gcd}(p'-q',q')=1$,
which is a contradiction.

Therefore, it follows from (\ref{rangeofalpha}) and the above claim that
\begin{equation}\label{smallest}
0<\{n\alpha^l\}(\alpha^k-1)<1
\end{equation}
for all $n\in\mathbb{Z}\setminus\{0\}.$
Then we conclude that 
\begin{align*}
    [[n\alpha^l]\alpha^k]+1 
    &=[n\alpha^{l+k} 
    -\{n\alpha^l\}\alpha^k]+1 \\
    &=[n\alpha^{l}+nm
    -\{n\alpha^l\}\alpha^k]+1 \tag{by (\ref{algfort})}\\
    &=[[n\alpha^{l}]-\{n\alpha^l\}(\alpha^k-1)] +1+nm\\
    &=[n\alpha^{l}]+nm \tag{by (\ref{smallest})}\\
    &= [n\alpha^{l+k}] \tag{by (\ref{algfort})}
\end{align*}
for all $n\in\mathbb{Z}\setminus\{0\}$.
\end{proof}

Now let us recall the well-known Kronecker-Weyl theorem for
the nontrivial direction of Theorem \ref{main}. The stated version below is much stronger than what we actually need; but it provides a deep insight into how the proof works.

\begin{theorem*}[Kronecker-Weyl]
    Let $\vec{\theta}=(\theta_1,\cdots,\theta_d)$ be a $d$-dimensional real vector, and write
    $$
    S_{\vec{\theta}}:=\left\{
    (\vec{a},b)\in\mathbb{Z}^d\times\mathbb{Q}:
    \vec{a}\cdot\vec{\theta}=b
    \right\}.
    $$
    Here we adopt the dot product 
    $$ \vec{a}\cdot\vec{\theta}=(a_1,\cdots,a_d)\cdot (\theta_1,\cdots,\theta_d):=
    a_1\theta_1+\cdots+a_d\theta_d.
    $$
Then \begin{enumerate}
    \item[(1)] Write $\pi_2:\mathbb{Z}^d\times\mathbb{Q}\to \mathbb{Q}$ to be the projection onto the second coordinate. Then all rational numbers in $\pi_2(S_{\vec{\theta}})$ have a least common denominator $q$.
    \item[(2)] Write 
    $$
    C:=\left\{\vec{x}\in\mathbb{R}^d: \vec{a}\cdot\vec{x}\equiv b\pmod{\mathbb{Z}}
    \text{ for all }(\vec{a},b)\in S_{\vec{\theta}}
    \right\}.
    $$
    Then the points in
    $$
    \bigcup\limits_{n\in\mathbb{Z}}\left\{\vec{x}\in\mathbb{R}^d:
    \vec{x}\equiv n\vec{\theta} \pmod{\mathbb{Z}^d}
    \right\}
    $$
    are equidistributed in the disjoint union of parallel affine subspaces
    $C\sqcup 2C\sqcup\cdots\sqcup qC$.
\end{enumerate}    
\end{theorem*}
\begin{proof}
    See \cite{w}, Satz 18.
\end{proof}

\begin{proof}[Proof of the ``if'' part of Theorem \ref{main}]
Since
\begin{align*}
    [[n\alpha^l]\alpha^k]+1 - [n\alpha^{l+k}]
    &= [n\alpha^{l+k}-\{n\alpha^{l}\}\alpha^{k}]+1-[n\alpha^{l+k}]\\
    &= [\{n\alpha^{l+k}\}-\{n\alpha^{l}\}\alpha^{k}]+1,
\end{align*}
    the identity (\ref{idforkl}) is equivalent to 
\begin{equation}\label{idforkl2}
    \{n\alpha^{l+k}\}-\{n\alpha^{l}\}\alpha^{k}\in [-1,0).
\end{equation} 
Define $x_n:=\{n\alpha^l\}$ and $y_n:=\{n\alpha^{l+k}\}$. Then (\ref{idforkl2}) is equivalent to the fact that the point $f(n):=(x_n,y_n)$ lies in the region 
\begin{equation}\label{regions}
   S:=\left\{(x,y)\in [0,1)^2: -1\leq y-\alpha^{k} x<0\right\}.
\end{equation}

Let us consider the orbit of discrete action
\begin{equation*}\label{orbit01}
    F:=\{f(n):n\in\mathbb{Z}\setminus\{0\}\}
\subseteq [0,1)^2.
\end{equation*}
By identifying $\mathbb{T}^2=\mathbb{R}^2/\mathbb{Z}^2$ with its fundamental domain $[0,1)^2$,
the subset $F\subseteq [0,1)^2$ corresponds to 
\begin{equation*}\label{orbittorus}
    F':=\{
(n\alpha^l,n\alpha^{l+k}) \pmod{\mathbb{Z}^2}: n\in\mathbb{Z}\setminus 
\{0\}
\}\subseteq \mathbb{T}^2.
\end{equation*}
In the following we will first study the properties of $F'$ and then come back to $F$. See Figure \ref{zzz} for
illustrations of the three cases we consider next.

If $\{1,\alpha^{l},\alpha^{l+k}\}$ is $\mathbb{Q}$-linearly independent, then it follows from Kronecker's theorem that $F'$ is dense in $\mathbb{T}^2$, namely, 
$F$ is dense in $[0,1)^2$. This leads to a contradiction since $F\subseteq S$.

Now suppose that $\{1,\alpha^{l},\alpha^{l+k}\}$ is $\mathbb{Q}$-linearly dependent. If $\alpha^l\in\mathbb{Q}$, then taking $n\in\mathbb{Z}\setminus \{0\}$ with $n\alpha^l\in\mathbb{Z}$ in (\ref{idforkl2}) also leads to a contradiction. If $\alpha^l\notin\mathbb{Q}$,
by $\mathbb{Q}$-linear dependence we have
\begin{equation}\label{qlindep}
\alpha^{l+k}=s\alpha^l+r
\end{equation}
for some $s,r\in\mathbb{Q}$.
We claim that the condition $F\subseteq S$ requires that 
\begin{equation*}
s=1\ \text{  and  }\ r\in\mathbb{Z}.
\end{equation*}

In fact, write $r=\frac{p}{q}$ where $p\in\mathbb{Z}$, $q\in\mathbb{N}$ and $\mathrm{gcd}(p,q)=1$. Then $\mathbb{Z}$ is divided into mod $q$ residue classes:
$$
\mathbb{Z}=\bigsqcup_{i=0}^{q-1}(q\mathbb{Z}+i).
$$
For each $i\in \{0,1,\cdots,q-1\}$, we compute the subsequence
of $F'\cup\{(0,0)\}$ given by
$$\{
((qm+i)\alpha^l,(qm+i)\alpha^{l+k})\pmod {\mathbb{Z}^2}:
m\in\mathbb{Z}
\}\subseteq \mathbb{T}^2$$
as follows:
    \begin{align*}
       \quad((qm+i)\alpha^l,(qm+i)\alpha^{l+k})\pmod {\mathbb{Z}^2}
       &=(qm+i) \cdot (\alpha^l,\alpha^{l+k}) \pmod{\mathbb{Z}^2}\\
       &=qm \cdot (\alpha^l,s\alpha^{l}+r) 
       +i\cdot(\alpha^l,\alpha^{l+k})\pmod{\mathbb{Z}^2}
       \tag{by (\ref{qlindep})}\\
       &=qm\alpha^l \cdot (1,s)
       +qm\cdot (0,r)
       +i\cdot(\alpha^l,\alpha^{l+k})\pmod{\mathbb{Z}^2}\\
       &\equiv qm\alpha^l\cdot (1,s)
       +i\cdot(\alpha^l,\alpha^{l+k})\pmod{\mathbb{Z}^2}.
    \end{align*}
Since $\alpha^l\notin\mathbb{Q}$, we see that
$\{
qm\alpha^l \pmod{\mathbb{Z}}: m\in\mathbb{Z}
\}$
is dense in $\mathbb{T}$.
With the condition $s\in\mathbb{Q}$,
this implies that 
$$\{
((qm+i)\alpha^l,(qm+i)\alpha^{l+k})\pmod {\mathbb{Z}^2}:
m\in\mathbb{Z}
\}$$
is dense in a translation of a one-dimensional subtorus in $\mathbb{T}^2$.
Therefore, we conclude that $F'\cup\{(0,0)\}$ (hence $F'$) is dense in a finite union of parallel subtori in $\mathbb{T}^2$, where the distances between neighboring subtori are the same.
This in turn means that $F$ is dense in a finite union of parallel segments with slope $s$ in $[0,1)^2$.

Let us take any one of these segments 
$$
L:=\{(x,y)\in [0,1)^2: y=sx+t\}
$$
for some $t\in\mathbb{R}$.
By the above paragraph, each of the segments 
\begin{align*}
    L'&:=\{(x,y)\in [0,1)^2: y=sx+t-1\},\\
    L''&:=\{(x,y)\in [0,1)^2: y=sx+t+1\},\\
    L'''&:=\{(x,y)\in [0,1)^2: y=s(x+1)+t\}
\end{align*}
must be a candidate for these segments unless it is empty.
We note that it follows directly from ($\ref{regions}$) that
$(L\cup L'\cup L''\cup L''')\setminus S$
will have nonempty interior in 
$L\cup L'\cup L''\cup L'''$, unless $s=1$ and $t=0$.
So the condition $F\subseteq S$, combined with the fact that $F\cap (L\cup L'\cup L''\cup L''')$ is dense in $L\cup L'\cup L''\cup L'''$,
requires that 
$$s=1\ \text{ and } \ t=0.$$
By arbitrariness we conclude that $F$ is contained and dense in 
$$\Delta:=\{(x,y)\in [0,1)^2:y=x\}.$$
In particular, we have $\{\alpha^{l}\}=\{\alpha^{l+k}\}$, namely,
    $\alpha^{l+k}-\alpha^{l}\in\mathbb{Z}$.
This verifies the claim.

\begin{figure}[H]\label{zzz}
    \centering
\begin{tikzpicture}
    \tikzmath{
        \length = 3; 
        \width = 0.5; 
        \slope = 1.618034; 
        \a = 4^(1/3); 
    }
    \clip (-\width,-\width) rectangle (\length+\width,\length+\width);
    \fill[black!20] (0,0)
        --($\length/\slope*(1,0)$)
        --($\length*(1,0)+{\length*(\slope-1)}*(0,1)$)
        --($\length*(1,1)$)
        --($\length/\slope*(1,0)+\length*(0,1)$)
        --cycle;
    \fill[black!5] (0,0)
        --($\length/\slope*(1,0)$)
        --($\length/\slope*(1,0)-\width*2/\slope*($(1,0)+\slope*(0,1)$)$)
        --($-\width*2/\slope*($(1,0)+\slope*(0,1)$)$)
        --cycle;
    \fill[black!5] ($\length*(1,1)$)
        --($\length*(1,0)+{\length*(\slope-1)}*(0,1)$)
        --($\length*(1,0)+{\length*(\slope-1)}*(0,1)+\width*2*($(1,0)+\slope*(0,1)$)$)
        --($\length*(1,1)+\width*2*(1,1)$)
        --($\length/\slope*(1,0)+\length*(0,1)+\width*2/\slope*($(1,0)+\slope*(0,1)$)$)
        --($\length/\slope*(1,0)+\length*(0,1)$)
        --cycle;

    \path[name path=solid] (0,\length)--(0,0)--(\length,0);
    \path[name path=dashed] (\length,0)--(\length,\length)--(0,\length);

    \foreach \n in {1,...,251}
    {
        \draw[fill,red] (${(\n*\a-int(\n*\a))*\length}*(1,0)+{(\n*\a^2-int(\n*\a^2))*\length}*(0,1)$) circle (0.01);
    }
    
    \draw[line width = 1pt] (0.01,\length)--(0,\length)--(0,0)--(\length,0)--(\length,0.01);
    \draw[line width = 1pt, dashed] (\length,0)--(\length,\length)--(0,\length);
    
    \draw[name path=straight,dashed] ($-\width*($(1,0)+\slope*(0,1)$)$)--(${(\length+\width)}*($(1,0)+\slope*(0,1)$)$);
    \path[name intersections={of=straight and dashed, by=Y1}];
    \node[below] (y1) at (Y1) {\tiny$y-\alpha^k x=0$}; 
    \draw[name path=straight] ($-\width*($(1,0)+\slope*(0,1)$)+($\length/\slope*(1,0)$)$)--(${(\length+\width)}*($(1,0)+\slope*(0,1)$)+($\length/\slope*(1,0)$)$);
    \path[name intersections={of=straight and dashed, by=Y2}];
    \node[below left] (y2) at (Y2) {\tiny$y-\alpha^k x=-1$}; 
    
    \node (S) at ($\length/\slope/2*(1,0)+0.5*($(1,0)+\slope*(0,1)$)$) {$S$};
\end{tikzpicture}
\begin{tikzpicture}
    \tikzmath{
        \length = 3; 
        \width = 0.5; 
        \slope = 1.618034; 
        \a = 1.914213562; 
    }
    \clip (-\width,-\width) rectangle (\length+\width,\length+\width);
    \fill[black!20] (0,0)
        --($\length/\slope*(1,0)$)
        --($\length*(1,0)+{\length*(\slope-1)}*(0,1)$)
        --($\length*(1,1)$)
        --($\length/\slope*(1,0)+\length*(0,1)$)
        --cycle;
    \fill[black!5] (0,0)
        --($\length/\slope*(1,0)$)
        --($\length/\slope*(1,0)-\width*2/\slope*($(1,0)+\slope*(0,1)$)$)
        --($-\width*2/\slope*($(1,0)+\slope*(0,1)$)$)
        --cycle;
    \fill[black!5] ($\length*(1,1)$)
        --($\length*(1,0)+{\length*(\slope-1)}*(0,1)$)
        --($\length*(1,0)+{\length*(\slope-1)}*(0,1)+\width*2*($(1,0)+\slope*(0,1)$)$)
        --($\length*(1,1)+\width*2*(1,1)$)
        --($\length/\slope*(1,0)+\length*(0,1)+\width*2/\slope*($(1,0)+\slope*(0,1)$)$)
        --($\length/\slope*(1,0)+\length*(0,1)$)
        --cycle;

    \path[name path=solid] (0,\length)--(0,0)--(\length,0);
    \path[name path=dashed] (\length,0)--(\length,\length)--(0,\length);

    \foreach \n in {1,...,150}
    {
        \draw[fill,red] (${(\n*\a-int(\n*\a))*\length}*(1,0)+{(\n*\a^2-int(\n*\a^2))*\length}*(0,1)$) circle (0.01);
    }
    
    \draw[line width = 1pt] (0.01,\length)--(0,\length)--(0,0)--(\length,0)--(\length,0.01);
    \draw[line width = 1pt, dashed] (\length,0)--(\length,\length)--(0,\length);
    
    \draw[name path=straight,dashed] ($-\width*($(1,0)+\slope*(0,1)$)$)--(${(\length+\width)}*($(1,0)+\slope*(0,1)$)$);
    \path[name intersections={of=straight and dashed, by=Y1}];
    \node[below] (y1) at (Y1) {\tiny$y-\alpha^k x=0$}; 
    \draw[name path=straight] ($-\width*($(1,0)+\slope*(0,1)$)+($\length/\slope*(1,0)$)$)--(${(\length+\width)}*($(1,0)+\slope*(0,1)$)+($\length/\slope*(1,0)$)$);
    \path[name intersections={of=straight and dashed, by=Y2}];
    \node[below left] (y2) at (Y2) {\tiny$y-\alpha^k x=-1$}; 

    \node (S) at ($\length/\slope/2*(1,0)+0.5*($(1,0)+\slope*(0,1)$)$) {$S$};
\end{tikzpicture}
\begin{tikzpicture}
    \tikzmath{
        \length = 3; 
        \width = 0.5; 
        \slope = 1.618034; 
        \a = 2.414213562; 
    }
    \clip (-\width,-\width) rectangle (\length+\width,\length+\width);
    \fill[black!20] (0,0)
        --($\length/\slope*(1,0)$)
        --($\length*(1,0)+{\length*(\slope-1)}*(0,1)$)
        --($\length*(1,1)$)
        --($\length/\slope*(1,0)+\length*(0,1)$)
        --cycle;
    \fill[black!5] (0,0)
        --($\length/\slope*(1,0)$)
        --($\length/\slope*(1,0)-\width*2/\slope*($(1,0)+\slope*(0,1)$)$)
        --($-\width*2/\slope*($(1,0)+\slope*(0,1)$)$)
        --cycle;
    \fill[black!5] ($\length*(1,1)$)
        --($\length*(1,0)+{\length*(\slope-1)}*(0,1)$)
        --($\length*(1,0)+{\length*(\slope-1)}*(0,1)+\width*2*($(1,0)+\slope*(0,1)$)$)
        --($\length*(1,1)+\width*2*(1,1)$)
        --($\length/\slope*(1,0)+\length*(0,1)+\width*2/\slope*($(1,0)+\slope*(0,1)$)$)
        --($\length/\slope*(1,0)+\length*(0,1)$)
        --cycle;

    \path[name path=solid] (0,\length)--(0,0)--(\length,0);
    \path[name path=dashed] (\length,0)--(\length,\length)--(0,\length);

    \foreach \n in {1,...,100}
    {
        \draw[fill,red] (${(\n*\a-int(\n*\a))*\length}*(1,0)+{(\n*\a^2-int(\n*\a^2))*\length}*(0,1)$) circle (0.01);
    }
    \draw[red,loosely dotted] ($\length*(0.5,0)$)--+($\length*(0,1)$);
    
    \draw[line width = 1pt] (0.01,\length)--(0,\length)--(0,0)--(\length,0)--(\length,0.01);
    \draw[line width = 1pt, dashed] (\length,0)--(\length,\length)--(0,\length);
    
    \draw[name path=straight,dashed] ($-\width*($(1,0)+\slope*(0,1)$)$)--(${(\length+\width)}*($(1,0)+\slope*(0,1)$)$);
    \path[name intersections={of=straight and dashed, by=Y1}];
    \node[below] (y1) at (Y1) {\tiny$y-\alpha^k x=0$}; 
    \draw[name path=straight] ($-\width*($(1,0)+\slope*(0,1)$)+($\length/\slope*(1,0)$)$)--(${(\length+\width)}*($(1,0)+\slope*(0,1)$)+($\length/\slope*(1,0)$)$);
    \draw[loosely dotted] ($\length/\slope*(1,0)$)--+($\length*(0,1)$);
    \path[name intersections={of=straight and dashed, by=Y2}];
    \node[below left] (y2) at (Y2) {\tiny$y-\alpha^k x=-1$}; 

    \node (S) at ($\length/\slope/2*(1,0)+0.5*($(1,0)+\slope*(0,1)$)$) {$S$};
\end{tikzpicture}
    \caption{The three possibilities to be excluded: the left image corresponds to the case when $ \{1,\alpha^l,\alpha^{k+l}\} $ is $ \mathbb{Q} $-linearly independent; the middle image corresponds to the case when the subtorus intersects the $ y $-axis except the origin; and the right image corresponds to the case when the slope is not equal to $ 1 $.}
\end{figure}
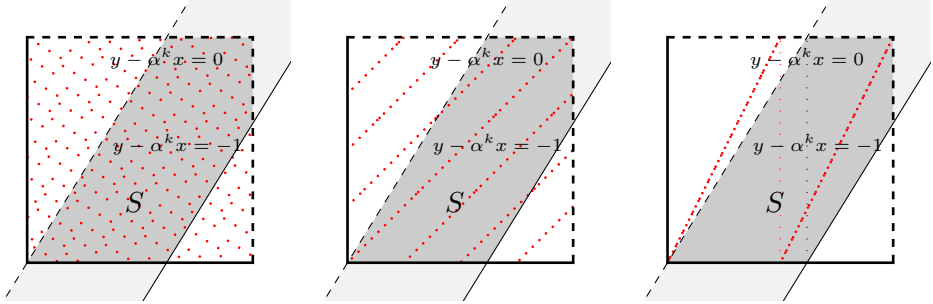

The remaining task is to estimate the range of $\alpha$. 
The above claim also requires that 
\begin{equation}\label{rangeofalpha2}
    \alpha^k-1\leq 1<\alpha^k,  
\end{equation}
otherwise $\Delta\setminus S$ would have nonempty interior in $\Delta$ and hence $F\not\subseteq S$.
(Here we also provide another elementary way to see (\ref{rangeofalpha2}):
On one hand, if $0<\alpha\leq 1$, then taking $n=1$ in $(\ref{idforkl2})$ leads to a contradiction, so 
$\alpha>1.$
On the other hand, if $\alpha^k>2$,
then by the irrationality of $\alpha^l$ we may choose some $n\in \mathbb{Z}\setminus\{0\}$ such that 
$\{n\alpha^l\}>\frac{2}{\alpha^k}$,
which leads to
$\{n\alpha^l\}\alpha^k>2>\{n\alpha^{l+k}\}+1$,
a contradiction to (\ref{idforkl2}).)

Moreover, if $\alpha=2^{\frac{1}{k}}$, then we may again
choose some $n\in \mathbb{Z}\setminus\{0\}$ such that $\{n\alpha^l\}<\frac{1}{2}$,
which leads to 
$\{n\alpha^{l+k}\}-\{n\alpha^{l}\}\alpha^{k}
=\{2n\alpha^l\}-2\{n\alpha^l\}=0$,
also a contradiction to (\ref{idforkl2}).
So we finally obtain
$$1<\alpha<2^{\frac{1}{k}},$$
and hence
$\alpha^{l+k}-\alpha^l\in [1,2^{\frac{l}{k}})$
as desired. 
\end{proof}

\section{Proof of Theorem \ref{gen1}} \label{proofofgen1}

The proof of Theorem \ref{gen1} is analogous to that of Theorem  \ref{main}.

\begin{proof}[Proof of Theorem \ref{gen1}]
The ``$\rm(i)\Rightarrow(iii)$'' part follows word by word as before once we notice that 
\[ 
0<\delta-\{n\alpha^l\}(\alpha^k-1)<1
\]
for all $n\in\mathbb{Z}$.

The ``$\rm(iii)\Rightarrow(i)$'' part is trivial.

    For the ``$\rm(ii)\Rightarrow(i)$'' part, the only difference to before is that the feasible region 
    becomes
  $$
         S:=\left\{(x,y)\in \mathbb{T}^2: 0\leq y-\alpha^{k} x+\delta< 1\right\}.
$$
With the same arguments, the requirement
    $F\subseteq S$
requires that $\alpha^{l+k}-\alpha^{l}\in\mathbb{Z}.$
Then the range of $\delta$ and $\alpha$ together
gives
$1<\alpha<2^{\frac{1}{k}}$.
\end{proof}

\section{Proof of Theorem  \ref{gen2}} \label{proofofgen2}


In this section, we first deduce Theorem \ref{gen2} and another interesting result from Proposition \ref{mostgen}, which is stated in a more general way. The proof of Proposition \ref{mostgen} will be split into two independent lemmas later.

\begin{proposition}\label{mostgen}
    Let $m$ be a positive integer, and let $\alpha$ and $\beta$ be two positive real numbers. Then the following statements are equivalent:
    \begin{enumerate}
        \item[\rm(i)] $\alpha>1$ is irrational,
    and $\beta-m\alpha\in\mathbb{Z}\cap [1,\alpha)$.
        \item[\rm(ii)] the identity
     \begin{equation}\label{m condition}
     [n\beta]-\left[ [n\alpha]\dfrac{\beta}{\alpha} \right]=[nm\alpha]+1-m[n\alpha]
    \end{equation}
    holds for all $n\in\mathbb{Z}\setminus\{0\}$. 
    \end{enumerate}
\end{proposition}

In particular, by taking $m=1$ in Proposition \ref{mostgen}, we have the following interesting result. 
\begin{corollary}
      A pair of positive real numbers $(\alpha,\beta)$ satisfies the identity
    \begin{equation*}
        [[n\alpha]\frac{\beta}{\alpha}]+1=[n\beta]
    \end{equation*}
    for all $n\in\mathbb{Z}\setminus\{0\}$, if and only if $\alpha>1$ is an irrational number, and
    \[\beta-\alpha\in \mathbb{Z}\cap [1,\alpha).\]
\end{corollary}

\begin{proof}[Proof of Theorem \ref{gen2} assuming Proposition \ref{mostgen}]
For the ``only if'' part, suppose that the condition (\ref{algforllm}) holds. Then it is easy to deduce that 
\begin{equation}\label{new rangeofalpha}
m<\alpha^k<m+1
\end{equation}
from the conditions $\alpha>0$ and $1\leq \alpha^l(\alpha^k-m)<(m+1)^{\frac{l}{k}}$.
Moreover, let us write
\begin{equation}\label{new algfort}
\alpha^{l}(\alpha^k-m)=M
\end{equation}
for some $M\in \mathbb{Z}\cap [1,(m+1)^{\frac{l}{k}})$ from (\ref{algforllm}).
We claim that 
\begin{equation}\label{new irrational}
\{\alpha^l,\alpha^k\}\cap\mathbb{Q}=\varnothing.
\end{equation}

In fact, suppose to the contrary that $\{\alpha^l,\alpha^k\}\cap\mathbb{Q}\neq\varnothing$.
Then it follows from (\ref{new algfort}) that 
$\{\alpha^l,\alpha^k\}\subseteq\mathbb{Q}$.
Let us write $\alpha^l=\frac{p}{q}$ and
    $\alpha^k=\frac{p'}{q'}$, where $\mathrm{gcd}(p,q)=\mathrm{gcd}(p',q')=1$.
Since $\alpha^k\in (m,m+1)$ cannot be an integer,
we have $\mathrm{gcd}(q,q')\neq 1$ by Lemma \ref{samedenom}.
However, it can be directly calculated from (\ref{new algfort}) that
$$
\frac{p}{q}=\frac{Mq'}{p'-mq'},
$$
which implies that $q\mid (p'-mq')$.
In particular, this means that 
$\mathrm{gcd}(q,q')\mid \mathrm{gcd}(p'-mq',q')=1$,
which is a contradiction.

Moreover, it follows from (\ref{new rangeofalpha}) and (\ref{new algfort}) that 
\begin{equation}\label{new algforllm3}
M=\alpha^l(\alpha^k-m)<\alpha^l.
\end{equation}
In summary, we conclude from (\ref{new rangeofalpha})-(\ref{new algforllm3}) and Proposition \ref{mostgen} that
$$
[n\alpha^{l+k}]-[ [n\alpha^l]\alpha^k ] =[nm\alpha^l]+1-m[n\alpha^l]
$$
for all $n\in\mathbb{Z}\setminus \{0\}$.

For the ``if'' part, suppose that the condition (\ref{idforllm}) holds. It follows from Proposition \ref{mostgen} that $$\alpha^{l+k}-m\alpha^l\in\mathbb{Z}\cap [1,\alpha^l).$$
In particular, we have $\alpha^k<m+1$ and hence
$\alpha^{l+k}-m\alpha^l\in\mathbb{Z}\cap [1,(m+1)^{\frac{l}{k}})$.

\end{proof}

The remaining part of this section is devoted to the proof of Proposition \ref{mostgen}. It is clear that Lemma \ref{converse theorem} provides the ``$\rm(i)\Rightarrow\rm(ii)$'' part of Proposition \ref{mostgen}, while Lemma \ref{new thm} provides the ``$\rm(ii)\Rightarrow\rm(i)$'' part of Proposition \ref{mostgen}.

\begin{lemma}
    \label{converse theorem}
    Let $m$ be a positive integer, $\alpha>1$ be an irrational number,
    and $\beta$ be a real number with 
    $\beta-m\alpha\in\mathbb{Z}\cap [1,\alpha)$. Then
\begin{equation*}
            [n\beta]-\left[ [n\alpha]\frac{\beta}{\alpha} \right] =[nm\alpha]+1-m[n\alpha]
\end{equation*}
for all $ n\in\mathbb{Z}\setminus\{0\}$.
\end{lemma}

\begin{proof}
Since $\beta-m\alpha\in\mathbb{Z}$, we have 
$\{n\beta\}=\{nm\alpha\}$,
namely,
\begin{equation}\label{m cond1}
    n\beta-[n\beta]=nm\alpha-[nm\alpha]
\end{equation}
for all $n\in\mathbb{Z}$.
Moreover, it follows from $\alpha\notin \mathbb{Q}$ and 
$\beta-m\alpha\in [1,\alpha)$ that 
\begin{equation}\label{m cond2}
    0<\{n\alpha\}\cdot\frac{\beta-m\alpha}{\alpha}<1
\end{equation}
for all $n\in\mathbb{Z}\setminus \{0\}$.
Then by direct calculation, we have
    \begin{align*}
        &\quad [n\beta]-\left[ [n\alpha]\frac{\beta}{\alpha} \right] \\
        &=n(\beta-m\alpha)+[nm\alpha]-\left[ [n\alpha]\frac{\beta-m\alpha}{\alpha} \right] -m[n\alpha]
        \tag{by (\ref{m cond1})}\\
        &=n(\beta-m\alpha)+[nm\alpha]-\left[ (n\alpha-\{n\alpha\})\frac{\beta-m\alpha}{\alpha} \right] -m[n\alpha]\\
        &=[nm\alpha]-\left[ -\{n\alpha\}\frac{\beta-m\alpha}{\alpha} \right] -m[n\alpha]\\
        &=[nm\alpha]+1-m[n\alpha]\tag{by (\ref{m cond2})}
    \end{align*}
    for all $n\in\mathbb{Z}\setminus \{0\}$.
\end{proof}

\begin{lemma}\label{new thm}
     Let $m$ be a positive integer, and let $\alpha$ and $\beta $ be two positive real numbers. Suppose that the identity
     \begin{equation*}
     [n\beta]-\left[ [n\alpha]\dfrac{\beta}{\alpha} \right]=[nm\alpha]+1-m[n\alpha]
    \end{equation*}
    holds for all $n\in\mathbb{Z}\setminus\{0\}$. 
    Then $\alpha>1$ is an irrational number, and
    \begin{equation}
    \beta-m\alpha\in\mathbb{Z}\cap [1,\alpha).
     \end{equation}
\end{lemma}

\begin{proof}

First of all, it is easy to see that the identity (\ref{m condition}) requires that $\alpha\notin\mathbb{Q}$. Otherwise taking $n\in\mathbb{Z}\setminus \{0\}$ with $n\alpha\in\mathbb{Z}$ in (\ref{m condition}) , the left-hand side is $0$ and the right-hand side is $1$, which is a contradiction.

Let us define
\begin{equation*}
    r(n):=[n\beta]-\left[ [n\alpha]\dfrac{\beta}{\alpha} \right]
\end{equation*}
for all $n\in\mathbb{Z}$. Then it follows from (\ref{m condition}) that
\begin{equation}\label{rn}
r(n)=[nm\alpha]+1-m[n\alpha]
\end{equation}
for all $n\in\mathbb{Z}\setminus \{0\}$.
We claim that 
\begin{equation}\label{m rnrange}
r(n)\in \{1,\cdots,m\}
\end{equation}
for all $n\in\mathbb{Z}\setminus\{0\}$, 
and moreover that,
for each $j\in \{1,2,\cdots,m\}$,
  \begin{equation}\label{m distribution}
        \lim_{N\to+\infty}\frac{|\{n\in\mathbb{Z}: 0<n\leq N,r(n)=j\}|}{N}=\frac{1}{m}.
  \end{equation} 

In fact, for any $n\in\mathbb{Z}\setminus \{0\}$, it is clear from (\ref{rn}) that $r(n)\in\mathbb{Z}$.
Since $m>0$, we have 
$[mx]\geq m[x]$ for all $x\in\mathbb{R}$,
and hence 
$$
r(n)=1+[m(n\alpha)]-m[n\alpha]\geq 1.
$$
Moreover, it follows from (\ref{rn}) that 
\begin{equation*}
    r(n)=-\{nm\alpha\}+1+m\{n\alpha\},
\end{equation*}
which implies that $r(n)<1+m$ and hence $r(n)\leq m$.
These verify (\ref{m rnrange}).

To verify (\ref{m distribution}), let us consider the orbit of discrete action
$$
F:=\{(\{n\alpha\},\{nm\alpha\}): n\in\mathbb{Z}\setminus \{0\}
\}\subseteq [0,1)^2.
$$
By identifying $\mathbb{T}^2=\mathbb{R}^2/\mathbb{Z}^2$ with the fundamental domain $[0,1)^2$, the subset $F\subseteq [0,1)^2$ corresponds to
$$
F':=\{(n\alpha,nm\alpha) \pmod{\mathbb{Z}^2}: n\in\mathbb{Z}\setminus \{0\}
\}\subseteq \mathbb{T}^2.
$$
Since $\alpha\notin\mathbb{Q}$, we see from Weyl's Equidistribution Theorem that 
$$\{n\alpha \pmod{\mathbb{Z}}:n\in\mathbb{Z}\}$$ is equidistributed in $\mathbb{T}$.
With the condition $m\in\mathbb{N}$, this implies that
$$
\{n\alpha\cdot (1,m)\pmod{\mathbb{Z}^2}:n\in\mathbb{Z}\}
$$
(and hence $F'$) is equidistributed in a one-dimensional subtorus of $\mathbb{T}^2$. This in turn means that $F$ is equidistributed in a finite union of parallel segments with slope $m$ in $[0,1)^2$, say
$\bigsqcup\limits_{j=1}^mL_j$, where 
$$
L_j:=\{(x,y)\in [0,1)^2: y=mx+1-j\}.
$$
It should be noted that for each $j\in\{1,2,\ldots,m\}$, we have
\begin{equation*}
(\{n\alpha\}, \{nm\alpha\})\in L_j
\Longleftrightarrow r(n)=j.
\end{equation*}
In particular, the equidistribution of $F$ in $\bigsqcup\limits_{j=1}^mL_j$ implies (\ref{m distribution}).

Now we turn to discussions around the two real numbers $\alpha$ and $\beta$.
To start with, we derive an alternative expression for $r(n)$ as follows: 
\begin{align*}
        r(n)&=[n\beta]-\left[ [n\alpha]\dfrac{\beta}{\alpha} \right]\\
        &=[n\beta]-\left[ (n\alpha-\{n\alpha\})\dfrac{\beta}{\alpha} \right]\\
        &=[n\beta]-\left[ n\beta-\{n\alpha\}\dfrac{\beta}{\alpha} \right]\\
        &=-\left[ \{n\beta\}-\{n\alpha\}\dfrac{\beta}{\alpha} \right].
    \end{align*}
Then for each $j\in \{1,\cdots,m\}$, the condition $r(n)=j$ is equivalent to the fact that the point $ (\{n\alpha\},\{n\beta\}) $ lies in the region
    \[
    S_j:=\left\{ (x,y)\in [0,1)^2: \left[ y-\frac{\beta}{\alpha}x \right] =-j \right\}.
    \]
Let us consider the orbit of discrete action 
$$
G:=\{
(\{n\alpha\},\{n\beta\}):n\in\mathbb{Z}\setminus \{0\}\}
\subseteq [0,1)^2,
$$
and the disjoint union of feasible regions
\begin{equation}\label{defs}
    S:=\bigsqcup\limits_{j=1}^mS_j=
\left\{ (x,y)\in [0,1)^2: -m\leq y-\frac{\beta}{\alpha}x<0 \right\}.
\end{equation}
Then it follows from (\ref{m rnrange}) that
\begin{equation}\label{gs}
    G\subseteq S.
\end{equation}
We claim that the condition (\ref{gs}) requires
\begin{equation}\label{intalphabeta}
    \beta\in\mathbb{Z}\alpha+\mathbb{Z}.
\end{equation}

In fact, by identifying $\mathbb{T}^2=\mathbb{R}^2/\mathbb{Z}^2$ with the fundamental domain $[0,1)^2$, the subset $G\subseteq [0,1)^2$ corresponds to
$$
G':=\{(n\alpha,n\beta) \pmod{\mathbb{Z}^2}: n\in\mathbb{Z}\setminus \{0\}
\}\subseteq \mathbb{T}^2.
$$
If $\{1,\alpha,\beta\}$ is $\mathbb{Q}$-linearly independent, then it follows from Kronecker's theorem that $G'$ is dense in $\mathbb{T}^2$, namely, 
$G$ is dense in $[0,1)^2$. This leads to a contradiction to (\ref{gs}).

Now suppose that $\{1,\alpha,\beta\}$ is $\mathbb{Q}$-linearly dependent.
Since $\alpha\notin\mathbb{Q}$, we have
$$
\beta=s\alpha+r
$$
for some $s,r\in\mathbb{Q}$.
A similar argument as in the proof of Theorem \ref{main} gives that
$G$ is equidistributed in a finite union of parallel segments with slope $s$ in $[0,1)^2$.

Let us take any one of these segments 
$$
I:=\{(x,y)\in [0,1)^2: y=sx+t\}
$$
for some $t\in\mathbb{R}$.
By the above paragraph, each of the segments 
\begin{align*}
    I_j'&:=\{(x,y)\in [0,1)^2: y=sx+t-j\},\quad j=1,2,\cdots,m\\
    I_j''&:=\{(x,y)\in [0,1)^2: y=sx+t+j\},\quad j=1,2,\cdots,m\\
    I'''&:=\{(x,y)\in [0,1)^2: y=s(x+1)+t\}
\end{align*}
must be a candidate for these segments unless it is empty.
We note that it follows directly from ($\ref{defs}$) that
$(\cup_{j=1}^m (I_j'\cup I_j'')\cup I''')\setminus S$
will have nonempty interior in 
$\cup_{j=1}^m (I_j'\cup I_j'')\cup I'''$, unless $s\in \mathbb{Z}\cap [1,m]$ and $ t\in \mathbb{Z}\cap [1-s,0]$.
So the condition $G\subseteq S$, combined with the fact that $G\cap (\cup_{j=1}^m (I_j'\cup I_j'')\cup I''')$ is dense in $\cup_{j=1}^m (I_j'\cup I_j'')\cup L'''$,
requires that 
\begin{equation}\label{rangest}
    s\in \mathbb{Z}\cap [1,m]\ \  \text{  and  } \ \ t\in \mathbb{Z}\cap [1-s,0].
\end{equation}
By arbitrariness we conclude that $G$ is contained and equidistributed in 
$\bigsqcup_{j=1}^sI_j$,
where
\begin{equation*}
    I_j:=\{
(x,y)\in [0,1)^2: y=sx+1-j
\}.
\end{equation*}
In particular, we have $\{\beta\}=s\{\alpha\}+1-j$ for some $j\in \{1,2,\cdots,s\}$, namely,
    \begin{equation}\label{intalphabetas}
        \beta-s\alpha=[\beta]-s[\alpha]+1-j\in\mathbb{Z}.
    \end{equation}
This verifies (\ref{intalphabeta}).

The above process also forces that
\begin{equation}\label{rangesalphabeta}
    s<\frac{\beta}{\alpha}\leq m+1,
\end{equation}
otherwise $(I_1\cup I_s)\setminus S$
will have nonempty interior in 
$I_1\cup I_s$ and hence $G\not\subseteq S$.

Moreover, we claim that 
\begin{equation}\label{s=m}
    s=m.
\end{equation}

In view of (\ref{rangest}), suppose to the contrary that $s<m$, namely $\frac{1}{s}>\frac{1}{m}$. Here on one hand, since $G$ is equidistributed in $\bigsqcup\limits_{j=1}^sI_j$, we have 
\begin{equation}\label{distributes}
\lim_{N\to+\infty}\frac{|\{n\in\mathbb{Z}: 0<n\leq N,
(\{n\alpha\},\{n\beta\})\in I_1
\}|}{N}=\frac{1}{s}.
\end{equation}
On the other hand, 
since 
$$
r(n)=1\Longleftrightarrow(\{n\alpha\},\{n\beta\})\in S_1,
$$
it follows from (\ref{m distribution}) that 
\begin{equation}\label{distributem}
\lim_{N\to+\infty}\frac{|\{n\in\mathbb{Z}: 0<n\leq N,
(\{n\alpha\},\{n\beta\})\in S_1
\}|}{N}=\frac{1}{m}.
\end{equation}
Thus we conclude from (\ref{distributes}) and (\ref{distributem}) that
$$
\mathrm{int}(I_1)\not\subseteq \overline{S_1},
$$
and that the intersection point of $\mathrm{int}(I_1)$ with $\partial{S_1}\cap \partial{S_2}$ should be 
$
\left(\frac{1}{m},\frac{s}{m}\right).
$
However, by inserting this coordinate into the segment 
$$\partial{S_1}\cap \partial{S_2}
=\left\{(x,y)\in [0,1]^2:y-\frac{\beta}{\alpha}x=-1\right\},$$
we obtain $$
\beta-s\alpha=m\alpha,
$$
which is a contradiction to (\ref{intalphabetas}) since $\alpha\notin\mathbb{Q}$. This verifies (\ref{s=m}).

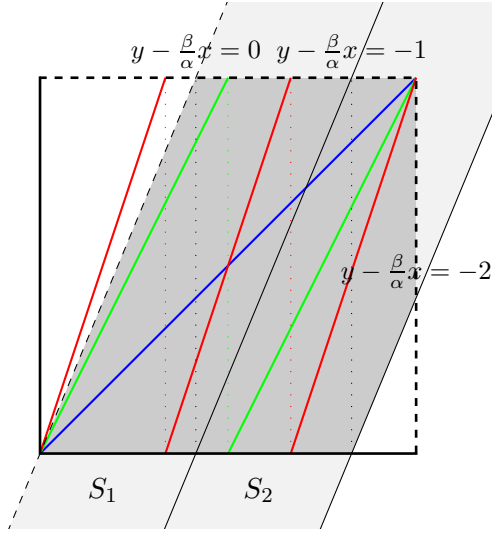
\begin{figure}[H]\label{fig2}
    \centering
\begin{tikzpicture}
    \tikzmath{
        \length = 5; 
        \width = 1; 
        \slope = 2.414213562; 
    }
    \clip (-\width,-\width) rectangle (\length+\width,\length+\width);
    \fill[black!20] (0,0)
        --($\length*2/\slope*(1,0)$)
        --($\length*(1,0)+{\length*(\slope-2)}*(0,1)$)
        --($\length*(1,1)$)
        --($\length/\slope*(1,0)+\length*(0,1)$)
        --cycle;
    \fill[black!5] (0,0)
        --($\length*2/\slope*(1,0)$)
        --($\length*2/\slope*(1,0)-\width*2/\slope*($(1,0)+\slope*(0,1)$)$)
        --($-\width*2/\slope*($(1,0)+\slope*(0,1)$)$)
        --cycle;
    \fill[black!5] ($\length*(1,1)$)
        --($\length*(1,0)+{\length*(\slope-2)}*(0,1)$)
        --($\length*(1,0)+{\length*(\slope-2)}*(0,1)+\width*2*($(1,0)+\slope*(0,1)$)$)
        --($\length*(1,1)+\width*2*(1,1)$)
        --($\length/\slope*(1,0)+\length*(0,1)+\width*2/\slope*($(1,0)+\slope*(0,1)$)$)
        --($\length/\slope*(1,0)+\length*(0,1)$)
        --cycle;
    \draw[thick, blue] (0,0)--(\length,\length);
    \draw[thick, green] (0,0)--(\length/2,\length);
    \draw[green, loosely dotted] (\length/2,0)--(\length/2,\length);
    \draw[thick, green] (\length/2,0)--(\length,\length);
    \draw[thick, red] (0,0)--(\length/3,\length);
    \draw[red, loosely dotted] (\length/3,0)--(\length/3,\length);
    \draw[thick, red] (\length/3,0)--(\length*2/3,\length);
    \draw[red, loosely dotted] (\length*2/3,0)--(\length*2/3,\length);
    \draw[thick, red] (\length*2/3,0)--(\length,\length);
    \draw[line width = 1pt] (0.01,\length)--(0,\length)--(0,0)--(\length,0)--(\length,0.01);
    \draw[line width = 1pt, dashed, name path=dashed] (\length,0)--(\length,\length)--(0,\length);
    \draw[name path=straight,dashed] ($-\width*($(1,0)+\slope*(0,1)$)$)--(${(\length+\width)}*($(1,0)+\slope*(0,1)$)$);
    \path [name intersections={of=straight and dashed, by=Y}];
    \node[above] (y) at (Y) {\small$y-\frac{\beta}{\alpha}x=0$}; 
    \foreach \n in {1,2}
    {
        \draw[name path=straight] ($-\width*($(1,0)+\slope*(0,1)$)+($\n*\length/\slope*(1,0)$)$)--(${(\length+\width)}*($(1,0)+\slope*(0,1)$)+($\n*\length/\slope*(1,0)$)$);
        \draw[loosely dotted] ($\n*\length/\slope*(1,0)$)--+($\length*(0,1)$);
        \path [name intersections={of=straight and dashed, by=Y}];
        \node[above] (y) at (Y) {\small$y-\frac{\beta}{\alpha}x=-\n$}; 
    }
    \node (S1) at ($\length/\slope/2*(1,0)-0.2*($(1,0)+\slope*(0,1)$)$) {$S_1$};
    \node (S2) at ($\length/\slope*3/2*(1,0)-0.2*($(1,0)+\slope*(0,1)$)$) {$S_2$};
\end{tikzpicture}
    \caption{An illustration of the proof when $m=2$. The orbit should be equidistributed in the parallel segments, and spend the same amount of time in different regions. Here the green segments are as desired; but the blue segment with slope $1$ is cut by the border between regions at some rational point, and the first of the red segments with slope $3$ doesn't lie in the region.}
\end{figure}

Finally, since $\alpha\notin\mathbb{Q}$, it follows from (\ref{intalphabetas}),(\ref{rangesalphabeta}) and (\ref{s=m}) that
$$
m=s<\frac{\beta}{\alpha}<m+1=s+1,
$$
namely, $\beta-m\alpha\in (0,\alpha)$.
See Figure \ref{fig2} for an illustration of the proof when m = 2.
\end{proof}


\bigskip

\section{Proof of Theorem  \ref{gen3}} \label{proofofgen3} 

The proof of Theorem \ref{gen3} is also analogous to that of Theorem  \ref{main}, with the help of Weyl's Equidistribution Theorem for polynomials from the Appendix.

\begin{proof}[Proof of Theorem \ref{gen3}]
The ``only if'' part follows from Theorem \ref{main}.

For the ``if'' part, let us write 
$$f(n):=(\{n\alpha^l\}, \{n\alpha^{l+k}\})$$ and 
$$S:=\left\{(x,y)\in \mathbb{T}^2: -1\leq y-\alpha^{k} x<0\right\}$$ 
as in the proof of Theorem \ref{main}, but redefine
$$F:=\{f(P(n)): n\in\mathbb{Z}\setminus\{
\text{roots of } P(X)\}
\}\subseteq [0,1)^2
$$
and 
$$F':=\{(P(n)\alpha^l,P(n)\alpha^{l+k})\pmod{\mathbb{Z}^2}:
n\in\mathbb{Z}\setminus\{
\text{roots of } P(X)\}
\}\subseteq \mathbb{T}^2.
$$
Then identity (\ref{idforpoly}) for all $n\in\mathbb{Z}\setminus\{
\text{roots of } P(X)\}$
is equivalent to
\begin{equation}\label{equiv}
    F\subseteq S.
\end{equation}

If $\{1,\alpha^{l},\alpha^{l+k}\}$ is $\mathbb{Q}$-linearly independent,
then by applying Corollary \ref{highdimpolyweyl} to 
$$
\widetilde{P}\colon \mathbb{Z}\to\mathbb{T}^2,\quad 
n\mapsto P(n)\cdot (\alpha^l,\alpha^{l+k}) \pmod{\mathbb{Z}^2},
$$
we see that $\{\widetilde{P}(n):n\in\mathbb{Z}\}$ (and hence $F'$) is dense in $\mathbb{T}^2$, namely, $F$ is dense in $[0,1)^2$. This leads to a contradiction to (\ref{equiv}).

Now suppose that $\{1,\alpha^{l},\alpha^{l+k}\}$ is $\mathbb{Q}$-linearly dependent. If $\alpha^l\in\mathbb{Q}$, then we write $\alpha^l=\frac{p}{q}$
where $p\in\mathbb{Z}$, $q\in\mathbb{N}$ and $\mathrm{gcd}(p,q)=1$.
Now $\mathbb{Z}$ is divided into mod $q$ residue classes:
$$\mathbb{Z}=\bigsqcup\limits_{i=0}^{q-1}(q\mathbb{Z}+i).$$
Write $P(n)=P_1(n)+a_0$ where $P_1(X)\in\mathbb{Z}[X]\setminus \{0\}$ with $P_1(0)=0$ and
$a_0\in\mathbb{Z}$.
For each $i\in\{0,1,\cdots,q-1\}$, we compute the subsequence
of $F'\cup\{(0,0)\}$ given by
$$\{\widetilde{P}(qm+i):m\in\mathbb{Z}\}\subseteq \mathbb{T}^2$$
as follows:
    \begin{align*}
       &\quad\widetilde{P}(qm+i)\\
       &=P(qm+i) \cdot (\alpha^l,\alpha^{l+k}) \pmod{\mathbb{Z}^2}\\
       &=P_1(qm+i) \alpha^l\cdot (1,0)
       +P_1(qm+i)\alpha^{l+k}\cdot (0,1)
       +a_0\cdot(\alpha^l,\alpha^{l+k})\pmod{\mathbb{Z}^2}\\
       &\equiv P_1(i) \alpha^l\cdot (1,0)
       +P_1(qm+i)\alpha^{l+k}\cdot (0,1)
       +a_0\cdot(\alpha^l,\alpha^{l+k})\pmod{\mathbb{Z}^2}.
    \end{align*}

Since $\alpha^l\in\mathbb{Q}$ and $\{\alpha^{l},\alpha^{k}\}\not\subseteq\mathbb{Q}$,
we have $\alpha^{l+k}\notin\mathbb{Q}$. Then it follows from Proposition \ref{onedimpolyweyl} that
$$\{P_1(qm+i) \alpha^{l+k} \pmod{\mathbb{Z}}:m\in\mathbb{Z}\}$$
is dense in $\mathbb{T}$.
This implies that 
$$\{\widetilde{P}(qm+i):m\in\mathbb{Z}\}$$
is dense in a translation of a vertical one-dimensional subtorus 
in $\mathbb{T}^2$.
Therefore, we conclude that
$F'\cup\{(0,0)\}$(and hence $F'$)
is dense in a finite union of vertical subtori in $\mathbb{T}^2$,
where the distances between neighboring subtori are the same.
This in turn means that $F$ is dense in a finite union of vertical segments in $[0,1)^2$, which is also a contradiction to (\ref{equiv}).

If $\alpha^l\notin\mathbb{Q}$, by $\mathbb{Q}$-linear independence we have
$$
    \alpha^{l+k}=s\alpha^{l}+r.
$$
for some $s,r\in\mathbb{Q}$. We claim that the condition $F\subseteq S$ requires that
$$
s=1\ \ \text{ and }\ \ r\in\mathbb{Z}.
$$

In fact, let us write $P(n)=P_1(n)+a_0$ where $P_1(X)\in\mathbb{Z}[X]\setminus \{0\}$ with $P_1(0)=0$,
$a_0\in\mathbb{Z}$,
and $r=\frac{p}{q}$ where $p\in\mathbb{Z}$, $q\in\mathbb{N}$, and $\mathrm{gcd}(p,q)=1$.
Then $\mathbb{Z}$ is divided into mod $q$ residue classes:
$$\mathbb{Z}=\bigsqcup\limits_{i=0}^{q-1}(q\mathbb{Z}+i).$$
For each $i\in\{0,1,\cdots,q-1\}$, we compute the subsequence
of $F'\cup\{(0,0)\}$ given by
$$\{\widetilde{P}(qm+i):m\in\mathbb{Z}\}\subseteq \mathbb{T}^2$$
as follows:
    \begin{align*}
       &\quad\widetilde{P}(qm+i)\\
       &=P(qm+i) \cdot (\alpha^l,\alpha^{l+k}) \pmod{\mathbb{Z}^2}\\
       &=P_1(qm+i) \cdot (\alpha^l,s\alpha^{l}+r) 
       +a_0\cdot(\alpha^l,\alpha^{l+k})\pmod{\mathbb{Z}^2}\\
       &=P_1(qm+i) \alpha^l\cdot (1,s)
       +P_1(qm+i)\cdot (0,r)
       +a_0\cdot(\alpha^l,\alpha^{l+k})\pmod{\mathbb{Z}^2}\\
       &\equiv P_1(qm+i) \alpha^l\cdot (1,s)
       +P_1(i)\cdot (0,r)
       +a_0\cdot(\alpha^l,\alpha^{l+k})\pmod{\mathbb{Z}^2}.
    \end{align*}

Since $\alpha^l\notin\mathbb{Q}$, we see from Proposition \ref{onedimpolyweyl} that
$$\{P_1(qm+i) \alpha^l \pmod{\mathbb{Z}}:m\in\mathbb{Z}\}$$
is dense in $\mathbb{T}$.
With the condition $s\in\mathbb{Q}$, this implies that 
$$\{\widetilde{P}(qm+i):m\in\mathbb{Z}\}$$
is dense in a translation of a one-dimensional subtorus 
in $\mathbb{T}^2$.
Therefore, we conclude that
$F'\cup\{(0,0)\}$(and hence $F'$)
is dense in a finite union of parallel subtori in $\mathbb{T}^2$.
The remaining arguments are the same as in Theorem \ref{main}.
\end{proof}

 \setcounter{secnumdepth}{0}

 \appendix
 \section{Appendix: Weyl's equidistribution criterion}
 \renewcommand{\thesection}{A}
 

The statement and proof of Weyl's classical equidistribution theorems
(Proposition \ref{Weylequi} and \ref{onedimpolyweyl})
can be found in \cite{w} or \cite{kn}. Here we formulate a version of Weyl's equidistribution theorem (Corollary \ref{highdimpolyweyl}) in the vector-valued polynomial case. 

\begin{proposition}[Weyl's Equidistribution Criterion]
\label{Weylequi}
Let $x\colon\mathbb{Z}\to \mathbb{T}^d$ be a function. Then $\{x(n):n\in\mathbb{Z}\}$ is equidistributed in $\mathbb{T}^d$ if and only if 
$$
 \lim\limits_{N\to +\infty} \dfrac{1}{N}\sum\limits_{|n|\leq N}
 e^{2\pi ik\cdot x(n)}=0
 $$   
for all $k\in \mathbb{Z}^d\setminus\{0\}$. Here we adopt the dot product
$$
(k_1,\cdots,k_d)\cdot (x_1,\cdots,x_d):=
k_1x_1+\cdots+k_dx_d \pmod{\mathbb{Z}}
$$
which maps $\mathbb{Z}^d\times \mathbb{T}^d$ to $\mathbb{T}$.
\end{proposition}

\begin{corollary}\label{reducedim}
    Let $x\colon\mathbb{Z}\to \mathbb{T}^d$ be a function. Then $\{x(n):n\in\mathbb{Z}\}$ is equidistributed in $\mathbb{T}^d$ if and only if, for each $k\in\mathbb{Z}^d\setminus\{0\}$,
    the sequence $\{k\cdot x(n):n\in\mathbb{Z}\}$ is equidistributed in $\mathbb{T}$.
\end{corollary}
\begin{proof}
    The result follows directly from Proposition \ref{Weylequi}.
\end{proof}

\begin{proposition}[Weyl's Equi-distribution Theorem for Polynomials]\label{onedimpolyweyl}   
    Let $l\geq 1$ be an integer, and let $P(X)=\alpha_lX^l+\cdots+\alpha_0$ be a polynomial of degree $l$ with $\alpha_0,\cdots,\alpha_l\in\mathbb{T}$. If $\alpha_l$ is irrational, then $\{P(n):n\in\mathbb{Z}\}$ is equidistributed in $\mathbb{T}$.
\end{proposition}

\begin{corollary}\label{highdimpolyweyl}
    Let $l\geq 1$ be an integer, and let $\widetilde{P}(X)=\alpha_lX^l+\cdots+\alpha_0$ be a polynomial of degree $l$ with $\alpha_0,\cdots,\alpha_l\in\mathbb{T}^d$. If $\alpha_l$ is irrational (i.e. for each $k\in\mathbb{Z}^d\setminus\{0\},
    k\cdot \alpha_l\neq 0$), then $\{\widetilde{P}(n):n\in\mathbb{Z}\}$ is equidistributed in $\mathbb{T}^d$.
\end{corollary}
\begin{proof}
    For each $k\in\mathbb{Z}^d\setminus \{0\}$,
    consider the polynomial 
    $$k\cdot \widetilde{P}(X)=(k\cdot\alpha_l)X^l+\cdots+(k\cdot\alpha_0)$$ with 
    $k\cdot \alpha_0,\cdots,k\cdot\alpha_l\in\mathbb{T}$.
    The condition requires that 
    $k\cdot\alpha_l\in\mathbb{T}\setminus\mathbb{Q}$.
    By Proposition \ref{onedimpolyweyl},
    we see that $n\mapsto k\cdot \widetilde{P}(n)$ is equidistributed in $\mathbb{T}$. Then it follows from Corollary \ref{reducedim} that $n\mapsto \widetilde{P}(n)$ is equidistributed in $\mathbb{T}^d$.
\end{proof}

\section*{Acknowledgment}

The authors would like to thank Nikolay Moshchevitin for his generous help, including but not limited to organizing the discussion group,
providing possible research directions,
and carefully reading several versions of this paper.
They also would like to thank Anton Shutov as well as other anonymous referees for their considerable comments and suggestions.


\begin{thebibliography}{1}

\bibitem{aa}
Boris Adamczewski and Jakub Konieczny.
\newblock Bracket words: a generalisation of {S}turmian words arising from generalised polynomials.
\newblock {\em Trans. Amer. Math. Soc.}, 376(7):4979--5044, 2023.

\bibitem{kn}
L.~Kuipers and H.~Niederreiter.
\newblock {\em Uniform distribution of sequences}.
\newblock Pure and Applied Mathematics. Wiley-Interscience [John Wiley \& Sons], New York-London-Sydney, 1974.

\bibitem{S21}
A.~V. Shutov.
\newblock Renormalizations of circle rotations.
\newblock {\em Chebyshevski\u i\ Sb.}, 5(4(12)):125--143, 2005.

\bibitem{w}
Hermann Weyl.
\newblock \"uber die {G}leichverteilung von {Z}ahlen mod. {E}ins.
\newblock {\em Math. Ann.}, 77(3):313--352, 1916.

\bibitem{av}
A.~A. Zhukova and A.~V. Shutov.
\newblock On two relations characterizing the golden ratio.
\newblock {\em Dal\cprime nevost. Mat. Zh.}, 21(2):194--202, 2021.

\bibitem{z}
V.~G. Zhuravlev.
\newblock One-dimensional {F}ibonacci tilings.
\newblock {\em Izv. Ross. Akad. Nauk Ser. Mat.}, 71(2):89--122, 2007.

\end{thebibliography}
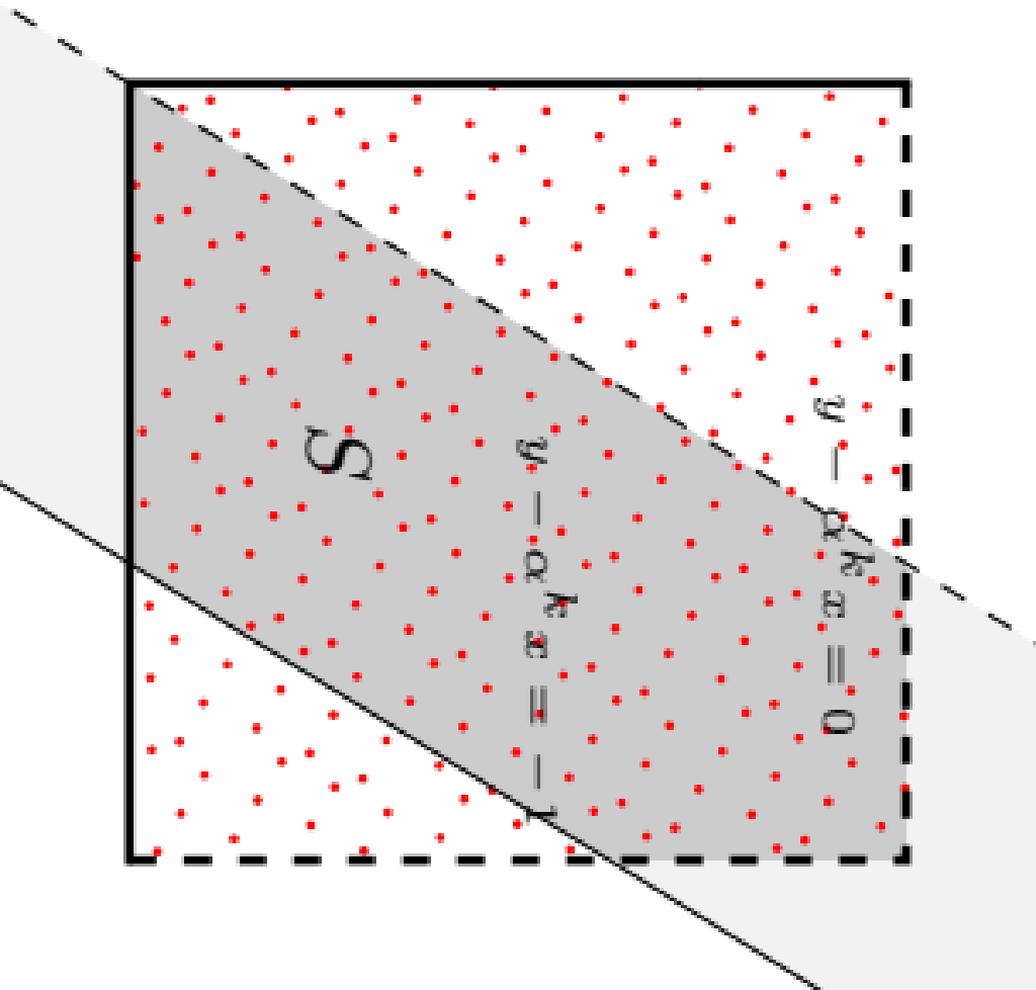
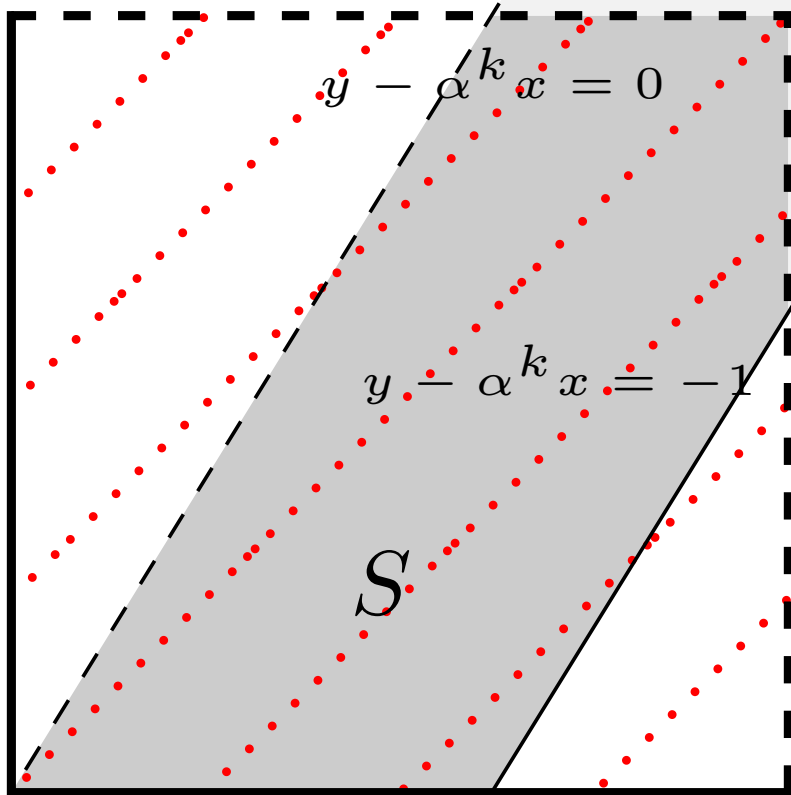
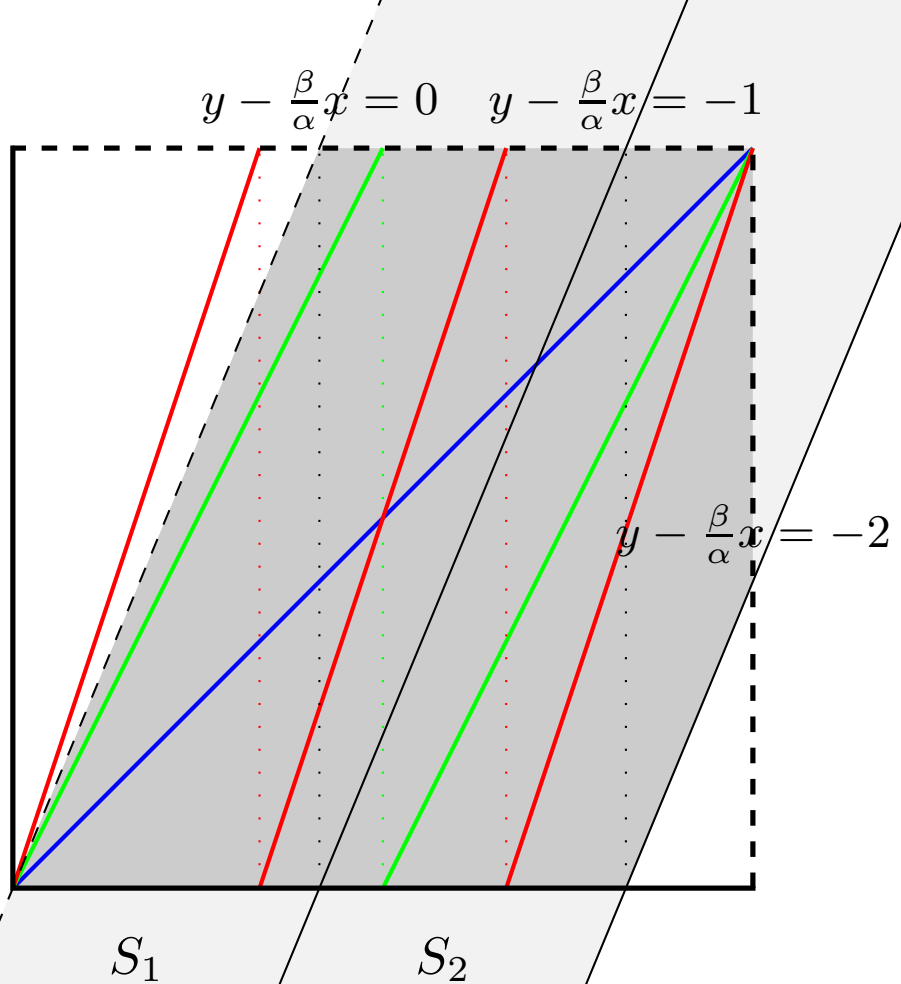
\end{document}